\documentclass[12pt]{article}
\usepackage{epic,eepic,epsf,epsfig}
\usepackage{amsmath,enumerate}
\usepackage{amssymb}
\usepackage{amsbsy}

\usepackage{amsfonts}
\newtheorem{theorem}{Theorem}[section]%
\newtheorem{lemma}[theorem]{Lemma}%
\newtheorem{cor}[theorem]{Corollary}%

 \def\Omega{\Omega}

\def\f{\noindent}

\newcommand{\qed}{\mbox{\raisebox{0.7ex}{\fbox{}}} \vspace{4truemm}}
\def\demo{\f {\bf Proof.}\hskip10pt}

\begin{document}

\baselineskip 16pt

\title{ \vspace{-1.2cm}
On a finite group in which every non-nilpotent maximal invariant subgroup of order divisible by $p$ is normal
\thanks{\scriptsize This work was supported in part by NSFC (12571024) and Shandong Provincial Natural Science Foundation, China (ZR2017MA022).
\newline
 \hspace*{0.5cm} \scriptsize $^{\ast\ast}$Corresponding
  author.\newline
 \hspace*{0.5cm} \scriptsize{E-mail addresses:} jtshi2005@pku.org.cn\,(J. Shi),\,\,lifh2026@s.ytu.edu.cn\,(F. Li),\,lina67@uzz.edu.cn\,(N. Li).}}

\author{Jiangtao Shi\,$^{a,\,\ast\ast}$,\,\,Fanghua Li$\,^a$,\,\,Na Li$\,^b$\\
\\
{\small{\em $^a$ School of Mathematics and Information Sciences,
Yantai University, Yantai 264005, P.R. China}}\\
{\small{\em $^b$ College of Mathematics and Statistics,  Zaozhuang University, Zaozhuang 277160, P.R. China}}}

\date{ }

\maketitle \vspace{-.8cm}

\begin{abstract}
Let $A$ and $G$ be finite groups such that $A$ acts coprimely on $G$ by automorphisms. For any fixed prime divisor $p$ of $|G|$, we prove that if
every non-nilpotent maximal $A$-invariant subgroup of $G$ of order divisible by $p$ is normal then $G$ has a Sylow tower and $G$ is $q$-nilpotent
or $q$-closed for each prime divisor $q$ of $|G|$.

\medskip \f {\bf Keywords:} non-nilpotent maximal invariant subgroup; normal; nilpotent; Sylow subgroup; Sylow tower
\end{abstract}

\section{Introduction}\label{s1}

In this paper all groups are considered to be finite. In {\rm\cite[Theorems 1.3 and 3.5]{lu2013}} Lu, Pang and Zhong showed that if every non-nilpotent
maximal subgroup of a group $G$ is normal then $G$ is $p$-nilpotent and $q$-closed for some primes $p$ and $q$. As a generalization, Shi {\rm\cite[Theorem 5]{shi2019}}
proved that a group $G$ in which every non-nilpotent maximal subgroup is normal has a Sylow tower. Furthermore, Shi, Li and Shen {\rm\cite[Theorem 1.7]{shi2023}} obtained
that such a group $G$ is either $q$-nilpotent or $q$-closed for any prime divisor $q$ of $|G|$. For any fixed prime divisor $p$ of $|G|$,
they {\rm\cite[Theorem 1.6]{shi2023}} proved that if every non-nilpotent maximal subgroup of $G$ of order divisible by $p$ is normal then there exists at most
one prime divisor $q\neq p$ of $|G|$ such that $G$ is neither $q$-nilpotent nor $q$-closed.

Combining the coprime action of groups, let $A$ and $G$ be groups such that $A$ acts coprimely on $G$ by automorphisms. Beltr\'{a}n and Shao {\rm\cite[Theorem A]{shao2023}} showed if
every non-nilpotent maximal $A$-invariant subgroup of $G$ is normal then $G$ is $p$-nilpotent for some prime $p$ dividing $|G|$. Li et al {\rm\cite[Theorem 1.5]{lu2023}}
had that such a group $G$ is a Sylow tower group. As an extension, Shi and Liu {\rm\cite[Theorem 1.5]{shiliu2025}
indicated that if every non-nilpotent maximal $A$-invariant subgroup of $G$ is normal then $G$ is $p$-nilpotent or $p$-closed for each prime divisor $p$ of $|G|$. Furthermore,
Shi and Xu {\rm\cite[Theorem 1.9]{shixu2025}} proved that every non-nilpotent maximal $A$-invariant subgroup of $G$ is normal if and only if
$G=(P_1\times\cdots\times P_{s-1})\times(P_s \rtimes(Q_1\times\cdots\times Q_t))$, where $P_1,\cdots,P_s$ are normal Sylow
subgroups of $G$ for $s\geq 1$, $Q_1,\cdots,Q_t$ are non-normal $A$-invariant Sylow subgroups of $G$
for $t\geq 1$, and there exists an $A$-invariant subgroup $E$ of $P_s$ such that $E$ is normal in $G$
and $E(Q_1\times\cdots\times Q_t)$ is a nilpotent maximal $A$-invariant subgroup of $P_s(Q_1\times\cdots\times Q_t)$.
Note that Beltr\'{a}n and Shao {\rm\cite[Theorem C]{shao2023}} showed that if every non-nilpotent maximal $A$-invariant subgroup of $G$ of order
divisible by $p$ is normal then $G$ is solvable, where $p$ is a prime divisor of $|G|$.

In this paper, as a further generalization and extension of above literature, we have the following result, the proof of which is given in Section~\ref{s3}.

\begin{theorem}\ \ \label{th00} Let $A$ and $G$ be groups such that $A$ acts coprimely on $G$ by automorphisms. If
every non-nilpotent maximal $A$-invariant subgroup of $G$ of order divisible by $p$ is normal for any fixed prime divisor $p$ of $|G|$. Then $G$ has a Sylow tower.
\end{theorem}

Furthermore, based on Theorem~\ref{th00}, we obtain the following result, the proof of which is given in Section~\ref{s4}.

\begin{theorem}\ \ \label{th01} Let $A$ and $G$ be groups such that $A$ acts coprimely on $G$ by automorphisms. If
every non-nilpotent maximal $A$-invariant subgroup of $G$ of order divisible by $p$ is normal for any fixed prime divisor $p$ of $|G|$. Then $G$
is either $q$-nilpotent or $q$-closed for each prime divisor $q$ of $|G|$.
\end{theorem}

The following corollary is a direct consequence of Theorem~\ref{th01}.

\begin{cor}\ \ \label{c1} If every non-nilpotent maximal subgroup of a group $G$ of order divisible by $p$ is normal for any fixed prime divisor $p$ of $|G|$. Then $G$
is $q$-nilpotent or $q$-closed for each prime divisor $q$ of $|G|$.
\end{cor}

Note that Theorem~\ref{th01} provides a complete generalization and improvement of {\rm\cite[Theorem 1.6]{shi2023}}.

\section{Several Necessary Lemmas}\label{s2}

\begin{lemma} {\rm\cite[Theorem B]{shao2019}}\ \ \label{l1} Let $G$ and $A$ be groups of coprime orders and assume that
$A$ acts on $G$ by automorphisms. If $G$ has a nilpotent maximal $A$-invariant subgroup of odd
order, then $G$ is solvable.
\end{lemma}

\begin{lemma} {\rm\cite[Lemma 2.3]{shao}}\ \  \label{lx} Suppose that a group $A$ acts coprimely on a group $G$. If every maximal $A$-invariant subgroup of $G$ is normal, then $G$ is nilpotent.
\end{lemma}

\begin{lemma} {\rm\cite[Lemma 9]{wie}}\ \ \label{l2} Let $H$ be a nilpotent Hall-subgroup of a group $G$ that is not
a Sylow subgroup of $G$. If for each prime divisor $p$ of $|H|$, assume $P\in{\rm Syl}_p(H)$ we always
have $N_G(P)=H$, then there exists a normal subgroup $K$ of $G$ such that $G=KH$ and $K\cap H=1$.
\end{lemma}

\section{Proof of Theorem~\ref{th00}}\label{s3}

\demo {\bf Part I}. We first provide a new proof to show that $G$ is solvable.

Let $G$ be a counterexample of the smallest order. For any $A$-invariant normal subgroup $N$ of $G$, it is easy to see that the hypothesis of the theorem
also holds for the quotient group $G/N$. It follows that $G$ has no non-trivial solvable $A$-invariant normal subgroups.

If every maximal $A$-invariant subgroup of $G$ of order divisible by $p$ is nilpotent. Then we can easily get that $G$ has normal Sylow subgroups arguing as in the proof
of {\rm\cite[Theorem 1.3]{shishanxu}} in which the solvability of the group is not involved. It implies that $G$ has non-trivial solvable $A$-invariant normal subgroups, a contradiction.

If every maximal $A$-invariant subgroup of $G$ of order divisible by $p$ is non-nilpotent. Then every maximal $A$-invariant subgroup of $G$ of order divisible by $p$ is
normal by the hypothesis. Let $P$ be an $A$-invariant Sylow $p$-subgroup of $G$. Since $P$ is not normal in $G$, $N_G(P)$ is a proper $A$-invariant subgroup of $G$. Then $G$ has a
maximal $A$-invariant subgroup $M$ such that $N_G(P)\leq M$. Note that $p\mid|M|$. By the assumption, $M$ is normal in $G$. Then $G=MN_G(P)=M$ by Frattini's argument, a contradiction.

It follows that $G$ not only has nilpotent maximal $A$-invariant subgroups of order divisible by $p$ but also has non-nilpotent maximal $A$-invariant subgroups of order divisible by $p$.

Let $H$ be a non-nilpotent maximal $A$-invariant subgroup of $G$ of order divisible by $p$ and $H_0$ a minimal $A$-invariant normal subgroup of $G$ such that $H_0\leq H$. Let $q$ be the largest prime divisor of
$|H_0|$. Since $H_0$ is non-solvable, $q>2$. Let $Q$ be an $A$-invariant Sylow $q$-subgroup of $G$. Then $Q\cap H_0$ is an $A$-invariant Sylow $q$-subgroup of $H_0$. Note that $Q\cap H_0$ is not normal
in $G$. It follows that $G$ has a maximal $A$-invariant subgroup $L$ such that $N_G(Q)\leq N_G(Q\cap H_0)\leq L$. In the following we discuss $L$ in three cases.

Case (1): Assume that $L$ is nilpotent.

Since $G$ has no non-trivial solvable $A$-invariant normal subgroups and $L$ is a nilpotent maximal $A$-invariant subgroup of $G$, $L$ must be a nilpotent Hall-subgroup of $G$. Moreover, let $R$ be any
$A$-invariant Sylow subgroup of $L$, one has $N_G(R)=L$. Note that
$q\mid|L|$ and $q$ is odd. Since $G$ is non-solvable, $L$ cannot be a Sylow subgroup of $G$ by Lemma\ref{l1}. Therefore, $G$ has a normal subgroup $E$ such that $G=E\rtimes L$ by Lemma~\ref{l2}.
Since $L$ is not only a maximal $A$-invariant subgroup of $G$ but also a Hall-subgroup of $G$, it follows that $E$ is not only a minimal normal subgroup of $G$ but also a Hall-subgroup of $G$.
Thus $E$ is characteristic in $G$, which implies that $E$ is also an $A$-invariant subgroup of $G$. Then $E$ is also a minimal $A$-invariant normal subgroup of $G$. Since $q\mid|L|$, one has
$q\nmid|E|$. It follows that $H_0$ and $E$ are two distinct minimal $A$-invariant normal subgroups of $G$, which implies that $H_0\cap E=1$. By the minimality of $G$, one has that both $G/H_0$ and $G/E$
are solvable. Then $G\cong G/(H_0\cap E)\lesssim(G/H_0)\times(G/E)$, which implies that $G$ is solvable, a contradiction.

Case (2): Assume that $L$ is non-nilpotent and $p\mid|L|$. By the hypothesis, $L$ is normal in $G$. It follows that $G=N_G(Q)L=L$ by Frattini's argument, a contradiction.

Case (3): Assume that $L$ is non-nilpotent and $p\nmid|L|$. Note that $G=N_G(Q\cap H_0)H_0=LH_0$ by Frattini's argument. One has $p\mid|H_0|$.
Let $P$ be an $A$-invariant Sylow $p$-subgroup of $G$. Then $P\cap H_0$ is an $A$-invariant Sylow $p$-subgroup of $H_0$. Note that $N_G(P)\leq N_G(P\cap H_0)<G$.
There exists a maximal $A$-invariant subgroup $K$ of $G$ such that $N_G(P)\leq N_G(P\cap H_0)\leq K$. If $K$ is non-nilpotent. Since $p\mid|K|$, one has $K\trianglelefteq G$ by the hypothesis.
Then $G=N_G(P)K=K$ by Frattini's argument, a contradiction. Thus $K$ is nilpotent.

Case $(i)$: When $p>2$ or $p=2$ but $K$ is not a Sylow subgroup of $G$, arguing as above Case (1), we can get a contradiction.

Case $(ii)$: When $p=2$ and $K=N_G(P)=P$ is an $A$-invariant Sylow 2-subgroup of $G$. By Frattini's argument, $G=H_0N_G(P\cap H_0)=H_0P$. Note that $G=H_0L$ and $p\nmid|L|$.
Thus one has $\frac{|H_0||P|}{|H_0\cap P|}=\frac{|H_0||L|}{|H_0\cap L|}$. Since $(\frac{|P|}{|H_0\cap P|},\frac{|L|}{|H_0\cap L|})=1$, it follows that $\frac{|P|}{|H_0\cap P|}=1$ and then $|P|=|H_0\cap P|$.
One has $P=H_0\cap P\leq H_0$. It follows that $G=H_0P=H_0$. However,
$H_0\leq H$ and $H$ is a maximal $A$-invariant subgroup of $G$, a contradiction.

All above analyses imply that the counterexample of the smallest order does not exist and so $G$ is solvable.

{\bf Part II}. We show that $G$ has normal Sylow subgroups.

Let $G$ be a counterexample of the smallest order and $N$ a minimal $A$-invariant normal subgroup of $G$. Since $N$ is a characteristic simple group of $G$ and $G$ is solvable, $N$ is an elementary abelian group
of prime-power order. Assume $N\cong {Z_q}^m$, where $q$ is a prime divisor of $|G|$ and $m\geq 1$. By the minimality of $G$, $G/N$ has a normal Sylow subgroup $RN/N$, where $R$ is an $A$-invariant Sylow $r$-subgroup
of $G$. Then $RN$ is normal in $G$. By our assumption, $R$ is not normal in $G$. Then $N_G(R)<G$. It follows that $G$ has a maximal $A$-invariant subgroup $M$ such that $N_G(R)\leq M$. By Frattini's argument,
one has $G=N_G(R)RN=N_G(R)N=MN$.

If $p\nmid|M|$. Then we must have $q=p$ and $N$ is a normal Sylow $p$-subgroup of $G$, a contradiction.

If $p\mid|M|$. We discuss $M$ in two cases.

Case $(i)$: Assume that $M$ is non-nilpotent. Then $M$ is normal in $G$ by the hypothesis. One has $G=N_G(R)M=M$ by Frattini's argument, a contradiction.

Case $(ii)$: Assume that $M$ is nilpotent. Let $Q_M$ be an $A$-invariant Sylow $q$-subgroup of $M$. Then $NQ_M$ is a normal Sylow $q$-subgroup of $G$, also a contradiction.

Hence the counterexample of the smallest order does not exist and then $G$ has normal Sylow subgroups.

{\bf Part III}. Final conclusion.

Let $Q_1$ be a normal Sylow $q_1$-subgroup of $G$. Considering the quotient group $G/Q_1$, the hypothesis of the theorem also holds for $G/Q_1$. Then $G/Q_1$ has a normal Sylow $q_2$-subgroup
$Q_2Q_1/Q_1$, where $Q_2$ is an $A$-invariant Sylow $q_2$-subgroup of $G$. An so on, we can get a normal subgroups series: $Q_1\vartriangleleft Q_1Q_2\vartriangleleft Q_1Q_2Q_3\vartriangleleft\cdots\vartriangleleft
Q_1Q_2Q_3\cdots Q_s=G$, where $Q_i$ is an $A$-invariant Sylow $q_i$-subgroup of $G$ for each $1\leq i\leq s$. It follows that $G$ has a Sylow tower. \hfill\qed

\section{Proof of Theorem~\ref{th01}}\label{s4}

\demo Let $P$ be an $A$-invariant Sylow $p$-subgroup of $G$. In the following we discuss $P$ in two cases.

Case (1): Assume that $P$ is normal in $G$. Let $R_1=P,\,R_2,\,\cdots,\,R_s$ be all normal Sylow subgroups of $G$, where $R_i$ is an $A$-invariant Sylow $r_i$-subgroup of $G$ for $1\leq i\leq s$.
Considering the $A$-invariant subgroup $U=P\times R_2\times\cdots\times R_s$. Since $G$ is solvable, $G$ has an $A$-invariant Hall-subgroup $M$ such that $G=U\rtimes M$.
Note that for any $A$-invariant Sylow subgroup $Q$ of $M$, one has that $Q$ is not normal in $G$ by the choice of above $R_i$ for $1\leq i\leq s$.

Claim that $M$ is nilpotent. Otherwise, assume that $M$ is non-nilpotent.
By Lemma~\ref{lx}, $M$ has an non-normal maximal $A$-invariant subgroup $H$. It is obvious that $UH$ is a maximal $A$-invariant subgroup of $G$.

Case $(i)$: Suppose that $UH$ is non-nilpotent. Then $UH$ is normal in $G$ by the hypothesis
since $p\mid|UH|$. It follows that $H=(U\cap M)H=UH\cap M$ is normal in $M$, a contradiction.

Case $(ii)$: Suppose that $UH$ is nilpotent. Then $UH=U\times H$.
It follows that $[U,H]=1$. For every $g\in M$, one has $[U,H]^g=[U^g,H^g]=
[U,H^g]=1$. It follows that $[U,H^M]=1$, where $H^M=\langle H^g|\,g\in M\rangle$ is a normal closure of $H$ in $M$. Since both $H$ and $M$ are $A$-invariant subgroups of $G$, one has that $H^M$ is also an $A$-invariant subgroup of $G$. It is easy to see that $H<H^M\leq M$ since $H$ is not normal in $M$. By the maximality of $H$, one has $H^M=M$.
Then $[U,M]=1$. It follows that $G=U\times M$. For any non-nilpotent maximal $A$-invariant $M_0$ of $M$, $U\times M_0$
is a non-nilpotent maximal $A$-invariant subgroup of $G$ of order divisible by $p$, one has that $U\times M_0$ is normal in $G$ by the hypothesis. It follows that
$M_0$ is normal $M$. Thus every non-nilpotent maximal $A$-invariant subgroup of $G$ is normal. By Theorem~\ref{th00}, $M$ has a normal Sylow subgroup $T$, where $T$ is an $A$-invariant Sylow subgroup of $G$.
Note that $[U,M]=1$. Then $[U,T]=1$. It follows that $U\leq N_G(T)$. Since $M\leq N_G(T)$ and
$G=U\times M$, one has that $T$ is normal in $G$, also a contradiction.

Therefore, $M$ is nilpotent. Hence $G$ is $q$-closed when $q=p,\,r_2,\,\cdots,\,r_s$ and $G$ is $q$-nilpotent when $q\mid|M|$.

Case (2): Assume that $P$ is not normal in $G$. Note that $G$ has normal Sylow subgroups. Let $Q_1,\,Q_2,\,\cdots,\,Q_t$ be all normal Sylow subgroups of $G$, where $Q_i$ is an $A$-invariant Sylow $q_i$-subgroup of $G$ for
$1\leq i\leq t$. Let $V=Q_1\times Q_2\times\cdots\times Q_t$. There exists an $A$-invariant Hall-subgroup $N$ of $G$ such that $G=V\rtimes N$. Note that $p\mid|N|$. For any non-nilpotent maximal $A$-invariant subgroup $N_0$ of $N$ of order divisible by $p$, $VN_0$ is a non-nilpotent maximal $A$-invariant subgroup of $G$ of order divisible by $p$. By the hypothesis, $VN_0$
is normal in $G$. It follows that $N_0$ is normal in $N$. Then $N$ has a normal Sylow subgroup $W$ by Theorem~\ref{th00}, where $W$ is an $A$-invariant Sylow subgroup of $G$. It follows that
$N\leq N_G(W)<G$. Let $K$ be a maximal $A$-invariant subgroup of $G$ such that $N_G(W)\leq K$. If $K$ is non-nilpotent. Since $p\mid|N|$, one has $p\mid|K|$. Then $K$ is normal in $G$ by the hypothesis.
It follows that $G=N_G(W)K=K$ by Frattini's argument, a contradiction. Thus $K$ is nilpotent, which implies that $N$ is nilpotent. Therefore, $G$ is $q$-closed when $q=q_1,\,q_2,\,\cdots,\,q_t$ and $G$ is $q$-nilpotent when $q\mid|N|$.\hfill\qed

\end{document}